\documentclass[12pt,a4paper]{amsart}
\usepackage[utf8]{inputenc}
\usepackage{tikz}
\usepackage[english]{babel}
\usepackage[T1]{fontenc}
\usepackage{amsmath}
\usepackage{amsfonts}
\usepackage{amssymb}
\usepackage{graphicx}
\usepackage[left=2.5cm,right=2.5cm,top=3cm,bottom=2.5cm]{geometry}
\usepackage{ulem}
\usepackage{color}
\newtheorem{defi}{Definition}[section]
\newtheorem{thm}{Theorem}[section]

\newtheorem{rem}{Remark}[section]

 \newtheorem{prop}{Proposition}[section]
\newtheorem{lem}{Lemma}[section]

\newtheorem{cor}{Corollary}[section]

\usepackage{color}

\begin{document}

 \title {On the Fundamental Eigenvalue Gap of
Sturm-Liouville Operators }

\address{ \newline
Team of Modeling and Scientific Computing, 
Department of Mathematics, 
Multidisciplinary Faculty of Nador, 
University of Mohammed First, Morocco}
\email{m.ahrami@ump.ac.ma}

\address{ \newline
Team of Modeling and Scientific Computing, 
Department of Mathematics, 
Multidisciplinary Faculty of Nador, 
University of Mohammed First, Morocco}
\email{z.elallali@ump.ma}

\address{School of Mathematics,
Georgia Institute of Technology,\
Atlanta GA 30332-0160, USA} 
\email{harrell@math.gatech.edu}




\author{Mohammed Ahrami, Zakaria El Allali, and Evans M. Harrell II}
 \maketitle


 \maketitle

\begin{abstract}
We use methods of direct optimization as in \cite{ElAHar} to find the minimizers of the fundamental gap of Sturm-Liouville operators on an interval,
{under the constraint that the potential is of single-well form and that the weight function is of single-barrier form, and under similar constraints expressed in terms of convexity.}

\end{abstract}



\section{Introduction}


The goal of this article is to find optimal estimates,
{under constraints on the form of the coefficient functions in \eqref{eP},
for}
the {\it fundamental eigenvalue gap}
$\Gamma := \lambda_2 - \lambda_1$ of the 
Sturm-Liouville equation
\begin{equation}\label{eP}
H(p,q) u :=  -\frac{d}{dx}
\left(p(x)\frac{du}{dx}\right)+V(x) u  = \lambda  w(x) u
\end{equation}
on a finite interval, with self-adjoint boundary conditions,
According to \cite{Wei}, \S 8.4, one may impose any separated homogeneous boundary conditions of the form
\begin{align*}
u(0) \cos \alpha - (pu')(0) \sin \alpha &= 0\\
u(\pi) \cos \beta - (pu')(\pi) \sin \beta &= 0,\nonumber
\end{align*}
{where $0\leq \alpha , \beta < \pi $, to make $H$ self-adjoint; the interval 
has been standardized as $[0, \pi]$ without loss of generality.  
To keep the exposition simple we restrict
to Dirichlet boundary conditions, i.e., $u(0) = u(\pi) = 0$, although 
our techniques also work with other self-adjoint boundary conditions, with
suitable changes.  Some further simplifications will be imposed below.
}

The quantity $\Gamma$
 is of interest as the ionization energy in quantum theory,
and sharp bounds for natural categories
of potential energies $V(x)$ that prevent its collapse, especially single-well and convex $V(x)$,
have been studied
since the 1980s; cf. \cite{ahrami,Ash1, Ash2,Ash3, Ash4, Cheng1, Bognar, Horvath1,Horvath2,Lavine} and references therein.
The problem of maximizing $\Gamma$ not only casts light on the physical problem of ionization, but is interesting as a mathematical problem in its own right.

Most prior work has assumed or emphasized the case where both $p(x)$ and the
weight $w(x)$ are held constant.
{Liouville transformations allow one to convert \eqref{eP} into equivalent equations
with different $p(x)$, $V(x)$, and $w(x)$, but in general only one of these three functions can
be eliminated.  In the first three sections of this article we standardize with $p(x) \equiv 1$, and recall how the Liouville transformation works in an appendix.
In summary, these remarks allow us to concentrate on the problem
\begin{equation}\label{ePnoP}
\left\{
\begin{array}{lr}
 -u''+V(x) u  = \lambda w(x)u, & x\in [0,\pi ]   \\
 u(0)=u(\pi)=0 & 
\end{array}
\right.
\end{equation}
in the following two sections.}

In \cite{Ash1}, for the problem \eqref{ePnoP} with Dirichlet conditions, Ashbaugh and Benguria proved that 
the optimal lower bound for $\Gamma$ for symmetric single-well potentials
is achieved if and only if $V$ is constant on $(0,\pi)$. In \cite{Lavine} Lavine considered the class of convex potentials on $[0, \pi]$ and proved,
with either Dirichlet or Neumann boundary conditions, that the constant potential function minimizes $\Gamma$.
Later Horv\'ath \cite{Horvath1} returned with Lavine's methods to the problem of single-well potentials,
but without symmetry assumptions, and again showed that the constant potential 
was optimal with some restrictions on the transition point, and in 2015 Yu and Yang \cite{YuYa}
extended Horv\'ath's result by allowing other transition points and both Dirichlet and Neumann conditions.
More recently El Allali and Harrell \cite{ElAHar} used direct optimization methods to prove
sharp lower bounds for $\Gamma$ 
with general
single-well potential $V(x)$, without any restriction on the transition point $a \in [0, \pi]$, 
and obtained similar results in the case where the potential is convex.  
El Allali and Harrell were furthermore able to analyze the case where $V = V_0 + V_1$, where $V_0$ is
a fixed background potential energy
and $V_1$ is assumed either single-well or convex.  In contrast to the earlier studies 
of single-well potentials, which restrict the transition point in one way or other,
the minimizing single-well
potentials they found
are in general step functions and {\it not} necessarily constant, unless extra conditions are imposed.
In the classic case where $p=1$ they recovered with different arguments the result of Lavine that
$\Gamma$ is uniquely minimized among convex $V$ by the constant, and in the case of
single-well potentials, with no restrictions on the position of the minimum,
they obtained a new, sharp bound, that  $\Gamma > 2.04575\dots$.
{Some further related articles are
Huang's discussion of the eigenvalue gap \cite{Huang2} and eigenvalue ratio \cite{Huang3} for the vibrating string with symmetric densities, i.e., allowing variable $p(x)$, and the works
of Ashbaugh and Benguria \cite{Ash3}, Huang and Law \cite{Huang4}, and Horv\`{a}th and Kiss \cite{Horvath2}, which include
other expressions related to the low-lying eigenvalues such as
eigenvalue ratios like $\frac{\lambda_n}{\lambda_1}$.}

\section{Simple properties of the fundamental gap $\Gamma$}

We shall use expressions such as $\lambda_k(V, w)$ and $\Gamma(V, w) = \lambda_2(V, w)- \lambda_1(V, w)$
to indicate the dependence on the fundamental gap on coefficients in {\eqref{ePnoP}} with respect to which we wish to optimize.  In this section 
we review and slightly extend some useful observations about $\Gamma$ that are familiar from previous sources such as \cite{Ash1,Lavine,ElAHar}.
Most importantly, there is an explicit formula for the first derivative with respect to perturbations of $V$ and $w$:

\begin{lem}\label{F-H}
Suppose that $V(.,t)$ and $w(.,t)$ are one-parameter families of real-valued, locally $L^1$ functions, with
$\inf V(x,\kappa) > - \infty$, $C \ge w(x,\kappa) \ge \frac 1 C$ for some $C > 0$, and
$\frac{\partial V}{\partial \kappa}(x,\kappa)$
and $\frac{\partial w}{\partial \kappa}(x,\kappa) \in L^1(0, \pi)$.
Then
\[
\frac{d\lambda_{n}(\kappa)}{d\kappa}=-\lambda_{n}\int_{0}^{\pi}\frac{\partial w}{\partial \kappa}(x,\kappa)u_{n}^{2}(x,\kappa)dx+\int_{0}^{\pi}\frac{\partial V}{\partial \kappa}(x,\kappa)u_{n}^{2}(x,\kappa)dx .
\]
\end{lem}

\begin{proof}  
Because $\frac{\partial V}{\partial \kappa}$ and $\frac{\partial w}{\partial \kappa}$ are relatively bounded perturbations, Kato's theory of 
analytic perturbations applies, and since the eigenvalues with separated homogeneous boundary conditions are simple, this 
justifies the use of a formal expansion to calculate the effect of the perturbation, \`a la Feynman-Hellmann:
Denoting $\dot{u}=\frac{\partial u}{\partial \kappa}$, differentiation of  {\eqref{ePnoP}} with respect to $\kappa$
gives
\[
\dot{u}''_{n}+(\dot{\lambda}_{n}w+\lambda_{n}\dot{w}-\dot{V})u_{n}
+(\lambda_{n}w-V)\dot{u}_{n}=0.
\]
We multiply the equation above by $u_{n}(x,\kappa)$ and integrate with respect to $x$ from $0$ to $\pi$.
This yields
\[
\int_{0}^{\pi}\dot{u}''_{n}u_{n}dx+\int_{0}^{\pi}(\lambda_{n}w-V)
\dot{u}_{n}u_{n}dx= -\int_{0}^{\pi}(\dot{\lambda}_{n}w
+\lambda_{n}\dot{w}-\dot{V})u_{n}^{2}dx.
\]
Observing that $u_n( \lambda_{n}w-V)=-u_{n}''  $,
$$ \int_{0}^{\pi}\dot{u}''_{n}u_{n}dx-\int_{0}^{\pi}\dot{u}_{n}u''_{n}dx= 
-\int_{0}^{\pi}(\dot{\lambda}_{n}w+\lambda_{n}\dot{w}-\dot{V})u_{n}^{2}dx . $$
Integrating by parts twice yields
$$[\dot{u}'_{n}u_{n}]_{0}^{\pi}-\int_{0}^{\pi}\dot{u}'_{n}u'_{n}dx
-[\dot{u}'_{n}u'_{n}]_{0}^{\pi}+\int_{0}^{\pi}\dot{u}'_{n}u'_{n}dx=
-\int_{0}^{\pi}(\dot{\lambda}_{n}w+\lambda_{n}\dot{w}-\dot{V})u_{n}^{2}dx ,   $$
with boundary conditions  $u_{n}(0,t)=u_{n}(\pi,t)=0$. This gives
$$  -\int_{0}^{\pi}(\dot{\lambda}_{n}w+\lambda_{n}\dot{w}-\dot{V})u_{n}^{2}dx =0,$$    
so that $$ \dot{\lambda}_{n} \int_{0}^{\pi}u_{n}^{2}wdx= -\lambda_{n}\int_{0}^{\pi}\dot{w}u_{n}^{2}dx+ \int_{0}^{\pi}\dot{V}u_{n}^{2}dx. $$
Noting that $ \displaystyle\int_{0}^{\pi}wu_{n}^{2}dx=1,$
$$  \dot{\lambda}_{n}   = -\lambda_{n}\int_{0}^{\pi}\dot{w}u_{n}^{2}dx+ \int_{0}^{\pi}\dot{V}u_{n}^{2}dx . $$       
\end{proof} 

We next adapt the monotonicity argument of \cite{Ash1,ElAHar}
to incorporate the weight:

\begin{lem}\label{lemma2-2}
Consider the problem  {\eqref{ePnoP}}  with the same assumptions on $V$ and $w$ as in Lemma \ref{F-H}.  Without loss of generality we standardize the first two normalized eigenfunctions
so that
$u_{1,2}(x) > 0$ for $0 < x < \epsilon$ for some $\epsilon$.  Then:
\begin{enumerate}
\item
$\dfrac{u_{2}}{u_{1}} $ is decreasing on $(0,\pi)$.

\item The equation $ \vert u_{1}(x)\vert=\vert u_{2}(x)\vert $ has
either one or two solutions on $(0,\pi)$.
 \item There exist two points $x_{-}$ and $x_{+}$
  $ 0 \leq x_{-}< x_{+} \leq \pi $, at least one of which is interior to $(0, \pi)$, such that 
$u_{1}^{2}(x)>u_{2}^{2}(x)$ on $(x_{-},x_{+})$ and $u_{1}^{2}(x)\le u_{2}^{2}(x)$ on $(x_{-},x_{+})^c$.
\item
The equation $ \lambda_1 \vert u_{1}^2(x)\vert=\lambda_2 \vert u_{2}^2(x)\vert $ has
either one or two
solutions on $(0,\pi)$.
 \item 
There exist two points $\widehat{x}_{-}$ and $\widehat{x}_{+}$
  $ 0 \leq \widehat{x}_{-}< \widehat{x}_{+} \leq \pi $, at least one of which is interior to $(0, \pi)$, such that 
$\lambda_1 u_{1}^{2}(x)> \lambda_2 u_{2}^{2}(x)$ on $(\widehat{x}_{-},\widehat{x}_{+})$ and $\lambda_1 u_{1}^{2}(x)\le \lambda_2 u_{2}^{2}(x)$ on $(\widehat{x}_{-},\widehat{x}_{+})^c$.
\end{enumerate}
\end{lem}

\begin{proof}
We first show that $\left(\dfrac{u_{2}}{u_{1}}\right)'<0$ for $0 < x < x_0$, where $u_2(x_0) = 0$ and hence that there can be at most 
one value $x_- \in (0, x_0)$ for which $u_1(x_-) = u_2(x_-)$.

The Wronskian is by definition
\[
W(x)=u_{1}(x)u_{2}'(x)- u_{2}(x)u_{1}'(x).
\]
Thus with the weight in {\eqref{ePnoP}},
\[
W'(x)=(\lambda_{1}-\lambda_{2})w(x)u_{1}(x) u_{2}(x),
\]
and for the the quotient $v(x):=\dfrac{u_{2}(x)}{u_{1}(x)}$,
\[
v'(x)= \frac{u_{1}(x)u_{2}'(x)- u_{2}(x)u_{1}'(x)}{u_{1}^{2}(x)} = 
\frac{W(x)}{u_{1}^{2}(x)}.
\]
Hence
\begin{equation}
\label{ratdecr}
v'(x)=\frac{1}{u_{1}^{2}(x)}\int_{0}^{x}(\lambda_{1}-\lambda_{2})w(t)u_{1}(t) u_{2}(t)dt<0,
\end{equation}
since $\lambda_{1}<\lambda_{2}$ and $ u_{1},u_{2}>0 $ on $(0,x_{0})$.

 Suppose that there exist distinct $\alpha_{1,2} \in(0,x_0)$ such that
\[
u_{2}(\alpha_{i}) = u_{1}(\alpha_{i}),~~~~i=1,2~ . 
\]
Then $v(\alpha_{1})=v(\alpha_{2})$.
By Rolle's theorem there exists $\xi \in(\alpha_{1},\alpha_{2})\subset (0,x_{0})$
such that$:$ $ v'(\xi)=0 $, but this contradicts \eqref{ratdecr}.

Since $u_2$ vanishes at a unique point $x_0$ and the same argument
can be carried out after the change of variables $x \to \pi - x$ and an adjustment of the sign of $u_2$, it follows that
$\left(\dfrac{u_{2}}{u_{1}}\right)$ is strictly monotonic on $(x_0, \pi)$ and that there is at most one 
value $x_+ \in (x_0, \pi)$ for which $u_1(x_+) = u_2(x_+)$.

At least one of the points $x_{\pm} \in (0, \pi)$, because if $x_- = 0$ and $x_+ = \pi$ then $u_1(x) > |u_2(x)|$ for all 
$x \in (0, \pi)$,  which would contradict $\|u_1\|_2 = \|u_2\|_2$.
\end{proof}

\section{Characterization of optimizers}
 
In this section, we determine the explicit form of the gap-minimizing potential and density function of problem {\eqref{ePnoP}}, 
closely following the strategy of \cite{ElAHar}.

\subsection{The class of single-well potentials and single-barrier densities}

\begin{defi}

Let $1<M \le \infty$. The function $V$ is called a 
{\rm single-well function} if $V$ is non-increasing on $[0,a]$ and non-decreasing on $[a,\pi]$,
for some $ a \in[0,\pi]$.
The point $ a $ is called a {\rm transition point} (with no assumption of uniqueness). 
The notation below will be used through this article:
$$ SW_{[0,\pi],M}=\lbrace V(x): 0\le V(x)\le M,\textrm{where } V 
\textrm{ is a single-well function on }[0,\pi] \rbrace .$$

\end{defi}

\begin{defi}
Let {$0 < N_< \le N_>  < \infty $}. The function $w$ is called a {\rm single-barrier {density}}
if $w$ is non-decreasing on $[0,b]$ and non-increasing on $[b,\pi]$
for some $b\in[0,\pi]$.
The following notation will be used through this article:
\[
SB_{[0,\pi],N<,N_>}=\left\{ w(x): N_< \leq w(x) \leq N_>, \textrm{where } w 
\textrm{ is a single-barrier density on } [0,\pi]\right\}.
\]
\end{defi}

\begin{defi}
Consider the Sturm Liouville problem
\begin{align}\label{s-l}
    {-u''+(V_0(x)+ V(x))u=\lambda w(x)u}&\\
    u(0)=u(\pi)=0&,\nonumber
\end{align}                                                               
where $V$ is a single-well function and $w$ is a single-barrier
{density. The background potential $V_0$  is assumed bounded and measurable.}
If there exist $V_{*} \in W_{[0,\pi],M}$ and 
$w_{*} \in  SB_{[0,\pi],N_<,N_>}$ such that
\[
\Gamma(V_{*},w_{*})= \inf{( \Gamma(V,w),V \in SW_{[0,\pi],M},
~w\in SB_{[0,\pi],N<,N_>})},
\]
then we call the function $V_{*}$ an {\rm optimal potential}  and the function  $w_{*}$ an {\rm optimal density} for problem \eqref{s-l}.
\end{defi}

{As in \cite{ElAHar}, we will use compactness of the sets $SW$ and $SB$, following from a theorem of Helly:}

\begin{prop}\label{compact}
For any sequence $f_{n}\in\Lambda$,  ($\Lambda=SW$ or $SB$ with any fixed positive $M, N_{<,>}$), there exist a subsequence $f_{n_{k}}$ and a function $f_\star$ such that  $f_{n_{k}}(x)\longrightarrow f_\star(x) \in\Lambda $ for a.e. $x$.
\end{prop}

\noindent
{For the proof,  see \cite{ElAHar},
Proposition $2.1$.}

\begin{cor}\label{existence}
There exist a potential $V_{*}\in SW_{[0,\pi],M} $ and a density $w_{*}\in SB_{[0,\pi],N<,N_>} $ that minimize $\Gamma[V,w]$.
\end{cor}

\begin{proof}  
According to Proposition \ref{compact},
since the gap $\Gamma[V,w]$ is positive, there exist minimizing sequences
$(V_{n})_{n\in \mathbb N}\subset\Lambda$ and  $(w_{n})_{n\in \mathbb N }\subset\Lambda$ such that
$$ \lim _{n\rightarrow +\infty}(\lambda_{2}(V_{n},w_{n})- \lambda_{1}(V_{n},w_{n})) = \inf   \lbrace \lambda_{2}(V,w)- \lambda_{1}(V,w);~~ V,w \in\Lambda  
\rbrace .$$
By Proposition \ref{existence} we may pass to subsequences in $\Lambda$
that converge pointwise a.e. to limits  $V_{*}\in\Lambda$ and $w_{*}\in\Lambda$:
$$\lim _{n\rightarrow +\infty} V_{n}=V_{*}  \quad \textsf{and}\quad
\lim _{n\rightarrow +\infty} w_{n}=w_{*} $$  By the dominated convergence theorem these sequences also converge in $L^{1}(0,\pi)$. Hence by continuity of $\Gamma$ with respect to relatively bounded perturbations,
 $$ \lim _{n\rightarrow +\infty}(\lambda_{2}(V_{n},w_{n})- \lambda_{1}(V_{n},w_{n})) = \inf   \lbrace \lambda_{2}(V,w)- \lambda_{1}(V,w);~~ V,w \in\Lambda  \rbrace 
=\lambda_{2}(V_{*},w_{*})- \lambda_{1}(V_{*},w_{*}) .$$
\end{proof}  

\begin{thm}\label{stepfn}
For any piecewise continuous, strictly positive weight function $w$, the optimal potential $V_*$ is a step function.  For any  
piecewise continuous potential function $V$, the optimal weight $w_*$ is a step function.  The same is true for jointly 
optimal $V_*$ and $w_*$.  In each case the optimizers have the following characterization:
\begin{enumerate}
\item
$V_*(x) = 0$ a.e. on a
connected component of $\{x: u_2^2(x) > u_1^2(x)\}$ and on the complement of that interval, $V_*(x) = \max(V_*)$ a.e.  
\item
$w_*(x) = N_>$ a.e. on a
connected component of $\{x: \lambda_2 u_2^2(x) > \lambda_1 u_1^2(x)\}$ and on the complement of that interval,  $w_*(x) = \min(w_*)$ a.e.  
\end{enumerate}

\end{thm}

\begin{proof}  
Let  $V_{*}\in SW_{[0,\pi],M}$ and 
$w_{*}\in SB_{[0,\pi],N<,N_>}$ be the minimizing potential and density guaranteed by the lemma.

We can characterize $V_{*}$ with the same argument as was used in \cite{ElAHar}, since Lemma \ref{lemma2-2}
allows for the possibility of variable weight.  We repeat it here to make the proof self-contained, and to point out a key difference.
By Lemma \ref{lemma2-2} there exist $x_{\pm}$: $0\leqslant x_{-}<x_{+}\leqslant \pi$, for which 
\begin{align*}
    &u_{2}^{2}(x)>u_{1}^{2}(x) \textrm{ on } (0,x_{-})\cup (x_{+},\pi)\quad\textrm{(one of these intervals may be vacuous)}\\
    &u_{1}^{2}(x)>u_{2}^{2}(x) \textrm{ on } (x_{-},x_{+}).
\end{align*}

We define a family of single-well potentials by perturbing $V_{*}$ so that
\[
V(x,\kappa)=\kappa V_{1}(x)+(1-\kappa)V_{*}(x) \quad t\in[0,1].
\]
Next, we work out the explicit form of $V_{*}$ in two cases, beginning with

\noindent
(i)
$x_{-} \le a < x_{+}$.
For definiteness we arrange by reflecting if necessary that
$0 < x_{-} \le a \le x_{+}$, $V_{\star}(x_+) \ge V_{*}(x_-)$, and $V_{*}$ is nondecreasing for $x \ge x_+$.  
 In this case we
proceed in two stages.  First, let
$$ V_{1}(x)=  \left\{
\begin{array}{c}
   V_{*}(x_{-}) \textrm{ on } (0,a)  
   \\
V_{*}(x_{+})\textrm{ on } (a,\pi).
\end{array}
\right.$$

We observe that $V_{1} \in SW_{[0,\pi],M}$ has been chosen so that
$V(x,\kappa) \in SW_{[0,\pi],M}$ 
and $V_{1}(x)$  has the opposite sign to $u_{2}^{2}(x,0)-u_{1}^{2}(x,0)$ a.e.

By the Feynman-Hellmann formula,  at $\kappa=0$,
\begin{equation}\label{FHforq}
\frac{d(\lambda_{2}(0)-\lambda_{1}(0))}{d \kappa}=\int_{0}^{\pi}[V_{1}(x)-V_{*}(x)][u_{2}^{2}(x,0)-u_{1}^{2}(x,0)]dx.
\end{equation}
Since
$V_{*}$ is a minimizer
and 
$$[V_{1}(x)-V_{*}(x)][u_{2}^{2}(x,0)-u_{1}^{2}(x,0)]dx\leqslant 0,$$
$$
0 \leqslant\frac{d(\lambda_{2}(0)-\lambda_{1}(0))}{dt} \leqslant 0,
$$
which implies that the integrand in \eqref{FHforq} equals 0 a.e.  Thus
$V_{1}(x)=V_{*}(x)$ a.e. on $[0,\pi]$.
I.e., $V_{*}(x) = V_{*}(x_{+}) \chi_{(a,\pi)}$ a.e.  But if $a > x_-$ the alternative choice
$$ V_{1}(x)=  \left\{
\begin{array}{c}
   0 \textrm{ on } (0,x_-)  
   \\
V_{*}(x_{+})\textrm{ on } (x_-,\pi)
\end{array}
\right.$$
is also valid, ensuring that $V(x,\kappa) \in SW_{[0,\pi],M}$ and that
$V_1(x)$ 
has the opposite sign to $u_{2}^{2}(x,0)-u_{1}^{2}(x,0)$ a.e.
We are thus led to the conclusion that $V_{*}(x) = V_{*}(x_{+}) \chi_{(x_,\pi)}$ a.e.  In particular, for $M>0$ the unique transition point
is $a = x_-$.

\noindent (ii)  Secondly, suppose that
$a <x_{-}$; the case when $x_{+}< a$, is similar. 
Let
$$ V_{1}(x)=  \left\{
\begin{array}{c}
V_{*}(a) \textrm{ on }(0,x_{-})
   \\
\,\,\,\,V_{*}(x_{+}) \textrm{ on } (x_{-},\pi).
\end{array}                                                                             \right.
$$
Then $V_{1} \in SW_{[0,\pi],M}$ with 

\begin{align*}
V_{1}(x)-V_{*}(x) &\ge 0 \textrm{ on } (0,x_{-})\cup(x_{+},\pi)\\
&\le 0 \textrm{ on } (x_{-},x_{+}).
\end{align*}

We note that
$V(x,t)\in SW_{[0,\pi],M}$,
and by the optimality of $V_{*}$,
$$
0 \leqslant \frac{d(\lambda_{2}(0)-\lambda_{1}(0))}{dt} \leqslant 0 .
$$
Hence, as before we conclude that
$V_{1}(x)=V_{*}(x)$ on $[0,\pi]$.  
Indeed, since we have concluded that $V_{*}(a) = V_{*}(x_-)$, we may as well redefine $a$ as $x_-$,
which reduces case (ii) to case (i).
In conclusion all optimal $V_{*}$ must be step functions with a unique jump coinciding with 
either $x_-$ or $x_+$.

In 
\cite{ElAHar}, where the weight was constant, it was possible to conclude that $V_{*}$ was of the form
$M \chi_I(x)$ for an interval $I$ by relying on the fact that adding a constant to the potential function does not change the gap $\Gamma$.
With a variable weight we are only able to conclude that in general $V_{*} = C \chi_I(x)$ for an undetermined constant $C\le M$.

It remains to characterize the optimal weight $w_{*}$ by a similar argument applied to $w$ for fixed $V$.
Just as before, we can use Lemma \ref{lemma2-2} to conclude that 
there exist $\widehat{x}_{\pm}$: $0\leqslant \widehat{x}_{-}<\widehat{x}_{+}\leqslant \pi$, satisfying 
\begin{align*}
    &\lambda_{2} u_{2}^{2}(x)>\lambda_{1} u_{1}^{2}(x) \textrm{ on } (0,\widehat{x}_{-})\cup (\widehat{x}_{+},\pi)\\
    &\lambda_{1}u_{1}^{2}(x)>\lambda_{2}u_{2}^{2}(x) \textrm{ on } (\widehat{x}_{-},\widehat{x}_{+}).
\end{align*}
Given an optimal $w_{*}$, we extend it to a family of weights by
$$ w(x,\kappa)=tw_{1}(x)+(1-\kappa)w_{*}(x) \quad t\in[0,1].
$$
There are again two cases, beginning with

\noindent
(i)
$\widehat{x}_{-} \le b <\widehat{x}_{+}$.
As before we
can arrange a convenient orientation by reflecting if necessary.  We thus posit that
$0 < \widehat{x}_{-} \le b \le \widehat{x}_{+}$, $w_{\star}(\widehat{x}_{+}) \le w_{*}( \widehat{x}_{-})$, and $w_{*}$ is nonincreasing for $x \ge \widehat{x}_{+}$.  
If we let
$$ 
w_{1}(x)=  \left\{
\begin{array}{c}
    w_*(\widehat{x}_{-}) \textrm{ on } (0,b) 
   \\
w_{*}(\widehat{x}_{+})\textrm{ on } (b,\pi),
\end{array}
\right.
$$
then $w(x,\kappa)\in SB_{[0,\pi],N<,N_>}$ 
and $w_{1}(x)$  has the opposite sign to $\lambda_1 u_{1}^{2}(x,0)-\lambda_2 u_{2}^{2}(x,0)$ a.e. 

Fixing $V$ and differentiating the eigenvalues at
$\kappa=0$ by the Feynman-Hellman formula,
$$ \frac{d(\lambda_{2}-\lambda_{1})}{d \kappa}=\int_{0}^{\pi}[w_{1}(x)-w_{*}(x)][\lambda_{1} u_{1}^{2}(x,0)-\lambda_{2}u_{2}^{2}(x,0)]dx. $$
By the optimality of $w_{*}$ and
since
$$[w_{1}(x)-w_{*}(x)][\lambda_{1}u_{1}^{2}(x,0)-\lambda_{2}u_{2}^{2}(x,0)]dx\leqslant 0,$$
$$
0 \leqslant\frac{d(\lambda_{2}(0)-\lambda_{1}(0))}{d \kappa} \leqslant 0
$$
implies that
$w_{1}(x)=w_{*}(x)$ a.e. on $[0,\pi]$.
\\
In conclusion the optimal $w_{*}$ must be a step function with at most one jump, located at $b$.
If $b > \widehat{x}_{-}$, then the alternative choice
$$ w_{1}(x)=  \left\{
\begin{array}{c}
   N_> \textrm{ on } (0,\widehat{x}_{-})  
   \\
w_{*}(\widehat{x}_{+})\textrm{ on } (\widehat{x}_{-},\pi)
\end{array}
\right.$$
still ensures that $w(x,\kappa) \in SW_{[0,\pi],M}$ and 
$w_1(x)$ 
and has the opposite sign to $\lambda_1 u_{1}^{2}(x,0)-\lambda_2 u_{2}^{2}(x,0)$ a.e.
We are thus led to the conclusion that  $w_{*}(x) = N_>$ a.e. for $x < \widehat{x}_{-}$ and  $w_{*}(x) = w_{*}(x_{+})$ a.e. for $x > \widehat{x}_{-}$.  In particular, for $N_> > N_<$ the unique transition point
is $b = \widehat{x}_{-}$.

\noindent (ii)  Secondly, suppose that
$b<\widehat{x}_{-}$; the case when $\widehat{x}_{+}<b$ is similar. 
Let
$$ w_{1}(x)=  \left\{
\begin{array}{c}
N_> \textrm{ on } (0,\widehat{x}_{-})   
   \\
\,\,\,\,w_{*}(\widehat{x}_{+}) \textrm{ on }( \widehat{x} _{-},\pi).
\end{array}                                                                                \right.$$
We have $w_{1}\in SB_{[0,\pi],N<,N_>}$ with 
\begin{align*}
w_{1}(x)-w_{*}(x) &\le 0 \textrm{ on } (0,\widehat{x}_{-})\cup(\widehat{x}_{+},\pi)\\
&\ge 0 \textrm{ on } (\widehat{x}_{-},\widehat{x}_{+}).
\end{align*}
Noting that
$w(x,t)\in SB_{[0,\pi],N<,N_>}$ and
using the optimality of $w_{*}$,
$$
0 \leqslant \frac{d(\lambda_{2}(0)-\lambda_{1}(0))}{d \kappa} \leqslant 0 .
$$
We conclude that
$w_{1}(x)=w_{*}(x)$ a.e. on $[0,\pi]$.  As with the characterization of $V$, this means that we can now suppose that $b =  \widehat{x} _{-}$ is a transition point, reducing this case to case (i).



\end{proof}                 

\begin{thm}\label{transcen}
The eigenvalues of
the Sturm-Liouville problem (\ref{ePnoP}) correspond to the real roots of the transcendental equation
   $$
\eta\tan(z(\pi-\widehat{x}_{-}))=-z\tan\left[\eta(\widehat{x}_{-}-x_{-})+\arctan\left(\frac{\eta}{t}\tan(tx_{-})\right)\right]~~~ if~~ \lambda > \frac{\max(V_{\star})}{\min(w_{\star})},
$$
where  $\eta:=\sqrt{\lambda N_>-\max(V_{\ast})}$, $z:=\sqrt{\lambda \min(w_{\ast})-\max(V_{\ast})}$, and $t:=\sqrt{\lambda N_>  }$.
\end{thm}
\begin{proof}  
By Theorem \ref{stepfn}, the optimal potential $V_{*}$ must be  of the form
 $$ V_{*}(x)=  \left\{
\begin{array}{c}
    0 \textrm{ on }(0,x_{-}) 
   \\
\max (V_{\ast})  \textrm{ on } (x_{-},\pi),
\end{array}                                                                                \right.$$
and the optimal density $w_{*}$ must be  of the form
$$ w_{*}(x)=  \left\{
\begin{array}{c}
 N_>     \textrm{ on } (0,\widehat{x}_{-})  
   \\
 \min (w_{\ast}) \textrm{ on } (\widehat{x}_{-},\pi).
                                                                                      \end{array}                                                                                \right.$$
The eigenfunctions are given by
$$ u(x)=  \left\{
\begin{array}{c}
\alpha_{1} \sin(tx)\textrm{ on }(0,x_{-})\hspace{5cm}
   \\ 
 \beta_{1}\sin(\eta(x-x_{-}))+\beta_{2}\cos(\eta(x-x_{-})) \textrm{ on }(x_{-},\widehat{x}_{-})
\\
  \alpha_{2} \sin(z(\pi-x))\textrm{ on }(\widehat{x}_{-},\pi).   
    \hspace{4cm}
           
\end{array}                                                                        \right.$$
Where
$ \eta=\sqrt{\lambda N_> -\max(V_{\ast})} $ , $z=\sqrt{\lambda \min(w_{\ast})-\max(V_{\ast})}$ , $t=\sqrt{\lambda N_>  }$ and $\alpha_1,\alpha_2,\beta_{1},\beta_{2}$ are real constants, with $\beta_{1} \neq 0$.

\noindent Continuity of $u$ at $x_{-}$ gives

$$
\beta_{2}=\alpha_{1}\sin(tx_{-})
$$
\noindent and continuity of $u'$ at $x_{-}$ gives

$$
\beta_{1}=\frac{t\alpha_{1}}{\eta}\cos(tx_{-}).
$$
Then 

$$
\frac{\beta_{2}}{\beta_{1}}=\frac{\eta}{t}\tan(tx_{-}).
$$

\noindent Continuity of $u$ at $\widehat{x}_{-}$ gives

$$\alpha_{1}\left( \frac{t}{\eta}\cos(tx_{-})   \sin(\eta(\widehat{x}_{-}-x_{-}))+\sin(tx_{-})\cos(\eta(\widehat{x}_{-}-x_{-})) \right)=\alpha_{2}\sin(z(\pi-\widehat{x}_{-})).   
$$
It follows that
 $$
\frac{t\alpha_{1}}{\eta}\cos(tx_{-})\sqrt{1+\left( \frac{\eta}{t}\tan(tx_{-})\right) ^2} \sin\left(  \eta(\widehat{x}_{-}-x_{-})+\arctan\left( \frac{\eta}{t}\tan(tx_{-})\right) \right)=\alpha_{2}\sin(z(\pi-\widehat{x}_{-})),  
$$
\noindent and by the continuity of $u'$ at $\widehat{x}_{-}$ ,

$$\alpha_{1}\eta\left( \frac{t}{\eta}\cos(tx_{-})   
\cos(\eta(\widehat{x}_{-}-x_{-}))-\sin(tx_{-})\sin
(\eta(\widehat{x}_{-}-x_{-})) \right)=-\alpha_{2}z\cos(z(\pi-\widehat{x}_{-})).   
$$
Therefore,

$$
\alpha_{1}t \cos(tx_{-})\sqrt{1+\left( \frac{\eta}{t}\tan(tx_{-})\right) ^2} \cos\left(  \eta(\widehat{x}_{-}-x_{-})+\arctan\left( \frac{\eta}{t}\tan(tx_{-})\right) \right)=-\alpha_{2}z\cos(z(\pi-\widehat{x}_{-})) .  
$$
Then
$$
\eta\tan(z(\pi-\widehat{x}_{-}))=-z\tan\left[\eta(\widehat{x}_{-}-x_{-})+\arctan\left(\frac{\eta}{t}\tan(tx_{-})\right)\right]. 
$$

\noindent This ends the proof.
\end{proof}

\section{Liouville transform of Sturm-Liouville operators}


In this section we apply Lavine's estimate on the fundamental gap to the Sturm-Liouville equation (\ref{ePnoP}),
 \begin{align*}
\begin{array}{lr}
 -u''+V(x) u  = \lambda w(x)u, & x\in [0,\pi ].  \\
\end{array}
\end{align*}

\noindent The eigenvalues of \eqref{ePnoP} coincide with these of the corresponding eigenvalue problem in Liouville normal form,

\begin{align}\label{GSM2}
\frac{d^{2}\eta}{d\xi^{2}}+(\lambda-\psi(\xi))\eta=0 ~~on~~[0,L],
\end{align}

\noindent where $\psi(\xi)$ is the Liouville potential defined by
\begin{align}\label{psi}
\psi(\xi)=\frac{w''}{4w^{2}}-\frac{5(w')^{2}}{16w^{3}}+\frac{V}{w},
\end{align}
with $L=\int_{0}^{\pi}\sqrt{w(t)}dt$. (For background on the Liouville transform we refer to \cite{Birk}).
In particular we have$:$ 
$$
\Gamma[V,w]=\Gamma[\psi].
$$

\begin{prop}\label{cvxLiou}
If the Liouville potential $\psi$ given by (\ref{psi}) of the Sturm-Liouville problem (\ref{GSM2}) is convex, then 
$$
\Gamma[V,w]\geqslant\frac{3\pi^{2}}{(\int_{0}^{\pi}\sqrt{w(t)}dt)^{2}},
$$
and equality is obtained if and only if $\psi$ is constant.

\end{prop}

\begin{proof}  
If the Liouville potential $\psi$ of the Sturm-Liouville problem (\ref{GSM2})
is convex,  then by \cite{Lavine}, we have $$\Gamma[V,w]\geqslant\frac{3\pi^{2}}{L^{2}}.$$
Because the interval $[0,\pi]$ is tranformed to $[0,L]$ with $L=\int_{0}^{\pi}\sqrt{w(t)}dt$. Thus
$$
\Gamma[V,w]\geqslant\frac{3\pi^{2}}{(\int_{0}^{\pi}\sqrt{w(t)}dt)^{2}}.
$$
\end{proof}  

\begin{rem}
To see that convexity of $V$ is not required for Proposition \ref{cvxLiou}, consider the example
$V(x)=-x^{2}$. Let $w(x)=x^{2}\neq \textrm{const}$, then {$\psi ^{''}(\xi)=\frac{3x^{2}-25}{x^{6}} \ge  0$ on  $[5,6]$.} 
On this interval the Liouville potential $\psi$ is convex and hence Proposition \ref{cvxLiou} is applicable, so
$$
\Gamma[V,w]\geqslant\frac{3\pi^{2}}{(\int_{5}^{6}\sqrt{w(t)}dt)^{2}},
$$
{\it i.e.}
$$  \Gamma[V,w]\geqslant 0.978803\dots.$$

\end{rem}

\begin{prop}
Consider the Sturm-Liouville problem \eqref{ePnoP} with positive
density function $ w \in C^{2}(0,\pi)$ and continuous convex potential $V$ on $[0,\pi]$.
If the fundamental gap satisfies
$$\Gamma[V,w] = 3\pi^{2}\left(\int_{0}^{\pi}\sqrt{w(t)}dt\right)^{-2}
$$
Then
$$
\frac{w^{(3)}}{4w^{2}}-\frac{w''w'}{2w^{3}}-\frac{10w^{(3)}}{16w^{3}}+\frac{15w'^{3}}{16w^{4}}+\frac{V'}{w}-\frac{Vw'}{w^{2}}=0.
$$

for all $  x \in [0,\pi] $.

\end{prop}

\begin{proof}



\noindent If the fundamental gap satisfies 
$$\Gamma[V,w]= 3\pi^{2}\left(\int_{0}^{\pi}\sqrt{w(t)}dt\right)^{-2}
$$
then the Liouville potential $\psi$  is constant, in which case $\psi'=0$. Therefore
$$\frac{w^{(3)}}{4w^{2}}-\frac{w''w'}{2w^{3}}-\frac{10w^{(3)}}{16w^{3}}+\frac{15w'^{3}}{16w^{4}}+\frac{V'}{w}-\frac{Vw'}{w^{2}}=0 .   $$
\end{proof}




\begin{prop}

The Liouville potential $\psi$ corresponding to  {\eqref{ePnoP}} is convex if 
\begin{align*}
\frac{w^{(4)}}{4w^{2}}&-\frac{w^{(3)}w'}{w^{3}}-\frac{w''^{2}}{2w^{3}}+\frac{3w''w'^{2}}{2w^{4}}-
\frac{10}{16}\left[\frac{w^{(4)}}{w^{3}}-\frac{3w^{(3)}w'}{w^{4}}\right]\\
&+
\frac{15}{16}\left[\frac{3w'^{2}w''}{w^{4}}-\frac{4w'^{4}}{w^{5}}\right]
+\frac{V''}{w}-\frac{V'w'}{w^{2}}
-
\frac{V'w'+Vw''}{w^{2}}+\frac{2w'^{2}V}{w^{3}}\ge 0
\end{align*}
on $[0,\pi]$.
\end{prop}

\paragraph{\it Proof} 

To analyze the  convexity of $\psi$, we will apply the nonnegativity criterion for
the second derivative. To differentiate the Liouville potential we make use of the
chain rule. In particular,
$$
\frac{d\psi}{d\xi}=\sqrt{g}\frac{d\psi}{dx}.
$$
If $g=1$ then $\psi$ is convex, so  $\displaystyle\frac{d^{2}\psi}{dx^{2}}\geqslant 0$.\\
As a consequence
$$
\frac{d\psi}{d\xi}=\frac{w^{(3)}}{4w^{2}}-\frac{w''w'}{2w^{3}}-\frac{10w^{(3)}}{16w^{3}}+\frac{15w'^{3}}{16w^{4}}+\frac{V'}{w}-\frac{Vw'}{w^{2}}.
$$
This yields that
\begin{align*}
\frac{d^{2}\psi}{d\xi^{2}}&=\frac{w^{(4)}}{4w^{2}}-\frac{w^{(3)}w'}{w^{3}}
-\frac{w''^{2}}{2w^{3}}+\frac{3w''w'^{2}}{2w^{4}}-
\frac{10}{16}\left[\frac{w^{(4)}}{w^{3}}-\frac{3w^{(3)}w'}{w^{4}}\right]
\\
&\quad +
\frac{15}{16}\left[\frac{3w'^{2}w''}{w^{4}}-\frac{4w'^{4}}{w^{5}}\right]
+\frac{V''}{w}-\frac{V'w'}{w^{2}}
-\frac{V'w'+Vw''}{w^{2}}+\frac{2w'^{2}V}{w^{3}}.
\end{align*}


\begin{thebibliography}{00}
\bibitem{ahrami}
\newblock { M. Ahrami and Z. El Allali,
\em Lower bounds on the fundamental spectral gap with Robin boundary conditions,}
\newblock  2021 UNC Greensboro PDE Conference. Electron. J. Diff. Eqns. Conf. 26 (2022), pp 1-11. 




\bibitem{Ash1}
\newblock {M. Ashbaugh and R. Benguria, 
\em Optimal lower bound for the gap between the first two eigenvalues of one-dimensional  Schr\"odinger operators with symmetric single-well potentials, }
\newblock Proc. Amer. Math. Soc. 105 (1989), 419--424.

\bibitem{Ash2}
\newblock {M. S. Ashbaugh and R. Svirsky, 
\em  Periodic potentials with minimal energy bands, }
\newblock Proc. Amer. Math. Soc., 114 (1992) 69--77.

\bibitem{Ash3}
\newblock {M. Ashbaugh and R. Benguria, 
\em Optimal bounds for ratios of eigenvalues of one-dimensional Schr\"odinger operators with Dirichlet boundary conditions and positive potentials,}
\newblock Comm. Math. Phys., 124 (1989), 403--415.

\bibitem{Ash4}
\newblock {M.~S. Ashbaugh, E.~M. Harrell II, and R. Svirsky,
\em On minimal and maximal eigenvalues gaps and their causes, }
\newblock Pac. J. Math. 147 (1991) 1--24.

\bibitem{Cheng1}
\newblock {Y. H. Cheng, S. Y. Kung, C. K. Law and W. C. Lian, 
\em The dual eigenvalue problems for the Sturm-Liouville system, }
\newblock Computers and Mathematics with Applications, 60 (2010) 2556--2563.



\bibitem{AsHa82}
\newblock {M.~S. Ashbaugh and E.~M. Harrell II,
\em Perturbation theory for shape resonances and large
barrier potentials,}
\newblock Commun. Math. Phys. 83 (1982) 151--170.



 
  \bibitem{Doob}
\newblock {J. L. Doob,
\em Measure theory,}
\newblock Graduate Texts in Mathematics (143), 1994 Springer.

\bibitem{ElAHar}
\newblock {Z. El Allali  and E.~M. Harrell II ,
{\em Optimal bounds on the fundamental spectral gap with single-well potentials}, }
\newblock Proc. Amer. Math. Soc., 150,(2022), 57--587.





\bibitem{Bognar}
\newblock {G. Bogn\`{a}r and O. Dosly,
 \em The ratio of eigenvalues of the Dirichlet eigenvalue problem for equations with one-dimensional p-Laplacian,}
\newblock Abstract and Applied Analyis, 2010 (2010).
 
 \bibitem{Kato}
\newblock {T. Kato,
\em Perturbation theory for linear operator,}
\newblock  1980 Springer-Verlag.

\bibitem{Huang1}
\newblock {M. J. Huang,
\em On the eigenvalue ratio for vibrating strings,}
\newblock Proc. Amer. Math. Soc., 127 (1999) 1805--1813.



\bibitem{Horvath1}
\newblock {M. Horv\'ath,
\em On the first two eigenvalues of Sturm-Liouville operators,}
\newblock Proc. Amer. Math. Soc. 131 4 (2002) 1215--1224.
 
\bibitem{Horvath2}
\newblock {M. Horv\'ath and M. Kiss,
\em A bound for ratios of eigenvalues of Schr\"{o}dinger operators with single-well potentials,}
\newblock Proc. Amer. Math. Soc., 134 (2006) 1425--1434.

\bibitem{Huang2}
\newblock {M. J. Huang,
\em The eigenvalue gap for vibrating strings with symmetric densities, }
\newblock Acta Math. Hungar., 117 (2007) 341--348.

\bibitem{Huang3}
\newblock {M. J. Huang,
\em A note on the eigenvalues ratio of vibrating strings, }
\newblock Acta Math. Hungar., 123 (2009) 265--271.

\bibitem{Huang4}
\newblock {M. J. Huang and C. K. Law,
\em Eigenvalue ratios for the regular Sturm-Liouville system,  }
\newblock Proc. Amer. Math. Soc., 124 (1996) 1427--1436.

 \bibitem{Lavine}
 \newblock {R. Lavine,
 \em The eigenvalue gap for one-dimensional convex potentials,}
 \newblock Proceedings of the American Mathematical Society 1994; 121, 
815--821. 
 
 
 
\bibitem{svir2}
\newblock {R. Svirsky,
\em Maximally resonant potentials subject to p-Norm constraints,}
\newblock Pac. J. Math. 129, 357--374  (1987).

\bibitem{YuYa}
\newblock {X. J. Yu and C. F. Yang,
\em The gap between the first two eigenvalues of Schr\"odinger operators with single-well potential,}
\newblock Appl. Math. Comp. 268 (2015) 275--283.

 \bibitem{Wei}
J. Weidmann,
\newblock {\em Linear operators in {H}ilbert spaces}, volume~68 of {\em
  Graduate Texts in Mathematics}.
\newblock Springer-Verlag, New York-Berlin, 1980.
\newblock Translated from the German by Joseph Sz{\"u}cs.

 \bibitem{Birk}
Garrett Birkhoff and Gian-Carlo Rota,
\newblock {\em  
Ordinary Differential Equations, 4th Edition.  
}
\newblock New York: Wiley, 1989.  See Chapter 10, section 9.



\bibitem{Zett}
A. Zettl
\newblock {\em Strum-Liouville theory},
\newblock American Mathematical Society, 2005.






 
\end{thebibliography}
\end{document}